\AddToHook{begindocument/before}{} 
\DeclareHookRule{begindocument}{hyperref}{before}{my-mdpi} 

\documentclass[preprints,article,accept,moreauthors,pdftex]{my-mdpi}

\firstpage{1} 
\pubvolume{13}
\issuenum{11}
\articlenumber{772}
\pubyear{2024}
\copyrightyear{2024}
\externaleditor{Academic Editor: Valery Y. Glizer}
\datereceived{23 May 2024} 
\daterevised{22 October 2024} 
\dateaccepted{5 November 2024} 
\datepublished{6 November 2024} 
\hreflink{https://doi.org/10.3390/\newline axioms13110772} 
\doinum{10.3390/axioms13110772}
\pdfoutput=1 

\Title{Optimal Control of Microcephaly Under Vertical Transmission \mbox{of Zika}}

\TitleCitation{Optimal Control of Microcephaly Under Vertical Transmission of Zika}

\Author{Dilara Yap\i \c{s}kan $^{1,2}$\orcidA{}, 
Cristiana J. Silva $^{1,3}$\orcidB{} and 
Delfim F. M. Torres $^{1,}$*\orcidC{}}

\AuthorNames{Dilara Yap\i \c{s}kan, Cristiana J. Silva and Delfim F. M. Torres}

\AuthorCitation{Yap\i \c{s}kan, D.; Silva, C.J.; Torres, D.F.M.}

\address{$^{1}$ \quad Center for Research and Development in Mathematics and Applications (CIDMA),
Department of Mathematics, University of Aveiro, 3810-193 Aveiro, Portugal; 
dilarayapiskan@ua.pt (D.Y.); cristiana.joao.silva@iscte-iul.pt (C.J.S.)\\
$^{2}$ \quad Department of Mathematics, Bal\i kesir University, Bal\i kesir 10145, Turkey\\
$^{3}$ \quad Department of Mathematics,  ISTA---School of Technology and Architecture, 
 Iscte---Instituto Universit\'{a}rio de Lisboa, 1649-026 Lisbon, Portugal
}

\corres{Correspondence: delfim@ua.pt}

\abstract{The Zika virus, known for its potential to induce neurological conditions 
such as microcephaly when transmitted vertically from infected mothers 
to infants, has sparked widespread concerns globally. 
Motivated by this, we propose an optimal control problem for the prevention 
of vertical Zika transmission. The novelty of this study lies in its consideration of 
time-dependent control functions, namely, insecticide spraying
and personal protective measures taken to safeguard pregnant women from 
infected mosquitoes. New results provide a way to minimize the number of 
infected pregnant women through the implementation of control strategies 
while simultaneously reducing both the associated costs of control measures
and the mosquito population, resulting in a decline 
in \mbox{microcephaly cases.}}

\keyword{\textls[-15]{optimal control; 
vertical transmission;
vector-borne diseases; 
Zika virus; 
microcephaly}}

\MSC{49M05; 92D30}

\begin{document}

\section{Introduction}

Zika virus is a mosquito-borne disease mainly transmitted to humans through the
bite of female mosquitoes \emph{{Aedes aegypti}}.  It presents a serious threat to 
public health due to its vertical transmission from
pregnant women to their babies, potentially resulting in heightened occurrences of 
neurological disorders like microcephaly \cite{WHO1,WHO2}. In February 2016, the 
World Health Organization (WHO) declared Zika-related microcephaly a Public Health 
Emergency of International Concern (PHEIC), confirming the causal link between the 
Zika virus and congenital malformations. The PHEIC was declared to have concluded 
by the WHO in November of the same year. Despite a global decline in cases of Zika 
virus disease from 2017 onwards, transmission continues at low levels in various 
countries across the Americas and other endemic regions \cite{WHO3}.

Mathematical modeling of infectious diseases is not only an important tool
in understanding the dynamics of the disease but also contributes to the
process of taking necessary measures to prevent disease transmission 
\cite{Martins,Cholera,Monkeypox}. Researchers have examined the vertical 
transmission of Zika as well as the development of microcephaly in newborn 
babies (for example, \cite{Agusto,Shah,Zika:Ndairou:etal}). Additionally, 
controlled mathematical models to identify the crucial characteristics 
causing the transmission may be more effective in predicting the future 
course of the epidemic and preventing transmission 
\cite{Mosquito,COVID-19,cancer,epidemicmodel}.
Therefore, various perspectives on optimal control strategies for Zika
transmission dynamics have been discussed in the \mbox{literature 
\cite{Agusto,Shah,Wang,Okyere,Ali}.}

{Our study proposes an optimal control problem for the prevention 
of vertical Zika transmission, which is a significant public health 
concern. Moreover, the proposed optimal control problem differs from all the others 
found in the literature on Zika transmission. Indeed, we consider 
the uncontrolled Zika transmission model in Brazil \cite{Zika:Ndairou:etal} 
and introduce time-dependent control functions representing personal protection 
and insecticide spraying.} The objective is to minimize the number of 
infected pregnant women through the implementation of control strategies 
while simultaneously reducing the associated costs of these control measures.

This paper is organized as follows: in Section~\ref{sec:zika:model}, we 
introduce the uncontrolled Zika model; in Section~\ref{sec:ocp}, we formulate the
 optimal control problem and derive the optimality system using the Pontryagin maximum 
principle; in Section~\ref{sec:numerical}, our focus is on conducting numerical simulations 
to showcase the impacts of optimal control strategies;
{and finally, we conclude with Section~\ref{sec:conclusions}.}


\section{Uncontrolled Zika Model}
\label{sec:zika:model}

In this section, we recall the main assumptions of the mathematical 
model for the spread of Zika virus as proposed in \cite{Zika:Ndairou:etal}.
The model considers women as the population under study. The total women population,
given by $N$, is subdivided into four mutually exclusive compartments,
according to disease status: susceptible pregnant women ($S$);
infected pregnant women ($I$); women who gave birth to babies without a
neurological disorder ($W$); and women who gave birth to babies with a
neurological disorder due to \mbox{microcephaly ($M$).}

As for the mosquitoes population, there are four state variables related to
the (female) mosquitoes: $A_{m}$, which corresponds to the
aquatic phase, which includes the egg, larva, and pupa stages; $S_{m}$, for
the mosquitoes that might contract the disease (susceptible); $E_{m}$,
for the mosquitoes that are infected but are not able to transmit the Zika
virus to humans (exposed); and $I_{m}$, for the mosquitoes capable of
transmitting the Zika virus to \mbox{humans (infected).}

\medskip 

The following assumptions are considered in our model:
\begin{enumerate}
\item[(A.1)] There is no immigration of infected humans;
\item[(A.2)] The total human populations $N$ is constant;
\item[(A.3)] \textls[25]{The coefficient of transmission 
of Zika virus is constant and does not vary \mbox{with seasons;}}
\item[(A.4)] After giving birth, pregnant women are no longer pregnant and
they leave the population under study at a rate $\mu_h$ equal to the rate of
humans birth;
\item[(A.5)] Death is neglected, as the period of pregnancy is much smaller
than the mean \mbox{humans lifespan;}
\item[(A.6)] There is no resistant phase for the mosquito due to its short
lifetime.
\end{enumerate}

Note that male mosquitoes are not considered in this study because they
do not bite humans and consequently they do not influence the dynamics of
the disease. The differential system that describes the model is composed 
of pregnant women and women who gave birth: 
\begin{equation}
\begin{cases}
\label{zikamodel1} \displaystyle{S^{\prime}(t) 
= \Lambda - \left(\phi B\beta_{mh} \frac{I_m(t)}{N} 
+ (1-\phi) \tau_1 + \mu_h\right) S(t)}, \\[2mm] 
\displaystyle{I^{\prime}(t) = \phi B\beta_{mh} \frac{I_m(t)}{N} S(t) 
- (\tau_2 + \mu_h) I(t)}, \\[2mm] 
\displaystyle{W^{\prime}(t) = (1-\phi) \tau_1 S(t) + (1-\psi)\tau_2 I(t) 
- \mu_h W(t)}, \\[2mm] 
\displaystyle{M^{\prime}(t) = \psi \tau_2 I(t) - \mu_h M(t)},
\end{cases}
\end{equation}
where $N = S(t) + I(t) + W(t) + M(t)$ 
is the total population (women), with $t \in [0, t_f]$. The parameter 
$\Lambda$ denotes the new pregnant women per week, $\phi$ stands for the
fraction of susceptible pregnant women that become infected, $B$ is the
average daily biting (per day), $\beta_{mh}$ represents the transmission
probability from infected mosquitoes $I_m$ (per bite), $\tau_{1}$ is the
rate at which susceptible pregnant women $S$ give birth (in weeks), 
$\tau_{2} $ is the rate at which infected pregnant women $I$ give birth (in
weeks), $\mu_{h}$ is the natural death rate for pregnant women, and $\psi$
denotes the fraction of infected pregnant women $I$ that give birth to babies
with a neurological disorder due to microcephaly. The above system 
\eqref{zikamodel1} is coupled with the dynamics of the mosquitoes: 
\begin{equation}
\begin{cases}
\label{zikamodel2} \displaystyle{A_m^{\prime}(t)
= \mu_b\left(1- \frac{A_m(t)}{K}\right) 
\left(S_m(t) + E_m(t) + I_m(t)\right) - \left(\mu_A + \eta_A\right)A_m(t)},\\[2mm] 
\displaystyle{S_m^{\prime}(t)= \eta_A A_m(t)- \Big( B\beta_{hm} \frac{I(t)}{N} 
+\mu_m \Big) S_m(t)}, \\[2mm] 
\displaystyle{E_m^{\prime}(t)= \Big( B\beta_{hm} \frac{I(t)}{N}\Big) S_m(t)
-(\eta_m + \mu_m) E_m(t)}, \\[2mm] 
\displaystyle{I_m^{\prime}(t)= \eta_m E_m(t) - \mu_m I_m(t)},
\end{cases}
\end{equation}
where parameter $\beta_{hm}$ represents the transmission probability from
infected humans $I_h$ (per bite), $\mu_b$ stands for the number of eggs at
each deposit per capita (per day), $\mu_A$ is the natural mortality rate of
larvae (per day), $\eta_A$ is the maturation rate from larvae to adult (per
day), $1/\eta_m$ represents the extrinsic incubation period (in days), 
$1/\mu_m$ denotes the average lifespan of adult mosquitoes (in days), and $K$
is the maximal capacity of larvae. See Table~\ref{statevar:parameters} for
the description of the state variables and parameters of the Zika model 
\eqref{zikamodel1}--\eqref{zikamodel2}. 
\begin{table}[H]
\caption{Variables and parameters of the Zika model 
\eqref{zikamodel1}--\eqref{zikamodel2}, 
{as given in \cite{Zika:Ndairou:etal}.}}
\label{statevar:parameters}{\ } 
\tabcolsep=.35cm
\begin{tabular}{ll}
\toprule
{\footnotesize {\textbf{Variable/Symbol}}} & {\footnotesize {\textbf{Description}}} \\ \midrule
{\footnotesize {$S(t)$}} & {\footnotesize {susceptible pregnant women}} \\ 
{\footnotesize {$I(t)$}} & {\footnotesize {infected pregnant women}} \\ 
{\footnotesize {$W(t)$}} & {\footnotesize {women who gave birth to babies
without a neurological disorder}} \\ 
{\footnotesize {$M(t)$}} & {\footnotesize {women who gave birth to babies with
a neurological disorder due to microcephaly}} \\ 
{\footnotesize {$A_{m}(t)$}} & {\footnotesize {mosquitoes in the aquatic phase}}\\ 
{\footnotesize {$S_{m}(t)$}} & {\footnotesize {susceptible mosquitoes}} \\ 
{\footnotesize {$E_{m}(t)$}} & {\footnotesize {exposed mosquitoes}} \\ 
{\footnotesize {$I_{m}(t)$}} & {\footnotesize {infected mosquitoes}} \\ \midrule
{\footnotesize {$\Lambda$}} & {\footnotesize {new pregnant women (per week)}} \\ 
{\footnotesize {$\phi$}} & {\footnotesize {fraction of $S$ that become infected}}
\\ 
{\footnotesize {$B$}} & {\footnotesize {average daily biting (per day)}} \\ 
{\footnotesize {$\beta_{mh}$}} & {\footnotesize {transmission probability from 
$I_m$ (per bite)}} \\ 
{\footnotesize {$\tau_{1}$}} & {\footnotesize {rate at which $S$ give birth (in
weeks)}} \\ 
{\footnotesize {$\tau_{2}$}} & {\footnotesize {rate at which $I$ give birth (in
weeks)}} \\ 
{\footnotesize {$\mu_{h}$}} & {\footnotesize {natural death rate}} \\ 
{\footnotesize {$\psi$}} & {\footnotesize {fraction of $I$ that gives birth to
babies with a neurological disorder}} \\ 
{\footnotesize {$\beta_{hm}$}} & {\footnotesize {transmission probability from 
$I_h$ (per bite)}} \\ 
{\footnotesize {$\mu_b$}} & {\footnotesize {number of eggs at each deposit per
capita (per day)}} \\ 
{\footnotesize {$\mu_A$}} & {\footnotesize {natural mortality rate of larvae
(per day)}} \\ 
{\footnotesize {$\eta_A$}} & {\footnotesize {maturation rate from larvae to
adult (per day)}} \\ 
{\footnotesize {$1/\eta_m$}} & {\footnotesize {extrinsic incubation period (in
days)}} \\ 
{\footnotesize {$1/\mu_m$}} & {\footnotesize {average lifespan of adult
mosquitoes (in days)}} \\ 
{\footnotesize {$K$}} & {\footnotesize {maximal capacity of larvae}} \\ \bottomrule
\end{tabular}
\end{table}

We consider system \eqref{zikamodel1}--\eqref{zikamodel2} with given initial
conditions:
\begin{equation*}
\begin{split}
&S(0) = S_0, \quad I(0) = I_0, \quad W(0) = W_0, \quad M(0) = M_0, \\
&A_m(0) = A_{m0}, \quad S_m(0) = S_{m0}, \quad E_m(0) = E_{m0}, \quad I_m(0)
= I_{m0},
\end{split}
\end{equation*}
with $\left( S_0, I_0, W_0, M_0, A_{m0}, S_{m0}, E_{m0}, I_{m0} \right) \in 
\mathbb{R}^{8}_{+}$. In what follows, we assume $\beta_{mh} = \beta_{hm}$.

The positivity and boundedness of solutions, as well as the existence and stability analysis 
of both disease-free and endemic equilibria, were studied in \cite{Zika:Ndairou:etal}. 
{Additionally, the basic reproduction number and its sensitivity was also 
analyzed in \cite{Zika:Ndairou:etal}}. 
Considering the data collected by the WHO 
between  4   February 2016 and  10 November  2016 in Brazil, the authors of 
\cite{Zika:Ndairou:etal} concluded that the parameters most sensitive to interventions 
are $B$ and $\beta_{mh}$. Hence, to mitigate the transmission of the Zika virus, 
it is imperative to implement control measures aimed at reducing the number 
of daily mosquito bites, $B$, and the transmission probability from the infected mosquitoes, 
$\beta_{mh}$. Furthermore, the fraction $\phi$, representing susceptible pregnant women ($S$) 
who contract the virus, exhibits a sensitivity index very close to $+1$. This underscores 
the critical importance of preventive measures aimed at safeguarding susceptible 
pregnant women from infection.

These conclusions encourage us to pursue the identification of optimal strategies 
for mitigating the transmission of the Zika virus. In the upcoming section, 
we address this by introducing an optimal control problem.


\section{Optimal Control Problem}
\label{sec:ocp} 

In this section, we formulate an optimal control problem for Zika transmission.  
Our objective is to minimize the number of infected pregnant women, reduce 
the mosquito population, and minimize the cost associated 
{with the implementation of the \mbox{control measures.}}

{According to the Centers for Disease Control and Prevention, 
the best way to prevent Zika is to be protected from mosquito bites 
\cite{CDC:protect:Zika}}. Moreover, everyone, including pregnant and breastfeeding women, 
should take steps to prevent mosquito bites \cite{CDC:protect:Zika}. When used as directed, 
EPA-registered insect repellents are proven safe and effective, even for pregnant 
and breastfeeding women. Taking this into account, we propose a controlled model 
by introducing in the Zika model \eqref{zikamodel1}--\eqref{zikamodel2} two control 
functions $u_1(\cdot)$ and $u_2(\cdot)$. The control $u_{1}$ represents protective 
clothing, insect repellent, and bed-nets to protect pregnant women from infected mosquitoes;
{while control ${u}_{2}$ refers to the insecticide spray applied to the mosquito population. 
The dynamical control system for Zika transmission that we
propose is then given by}
\begin{equation}
\begin{cases}
S^{\prime}(t)=\Lambda -\left( 1-u_{1}(t)\right) {\phi B\beta _{mh}
\frac{I_{m}}{N}S(t)-(1-\phi )\tau _{1}S(t)-\mu _{h}S(t)}, \\ 
I^{\prime}(t)=\left( 1-u_{1}(t)\right) \phi B\beta _{mh}\frac{I_{m}(t)}{N}S(t)-(\tau
_{2}+\mu _{h})I(t), \\ 
W^{\prime}(t)=(1-\phi )\tau _{1}S(t)+(1-\psi )\tau _{2}I(t)-\mu _{h}W(t), \\ 
M^{\prime}(t)=\psi \tau _{2}I-\mu _{h}M(t), \\ 
A_{m}^{\prime}(t)=\mu _{b}\left( 1-\frac{A_{m}(t)}{K}\right) \left(
S_{m}(t)+E_{m}(t)+I_{m}(t)\right) -\left( \mu _{A}+\eta _{A}\right) A_{m}(t), \\ 
S_{m}^{\prime}(t)=\eta _{A}A_{m}(t)-B\beta _{hm}\frac{I(t)}{N}S_{m}(t)+\mu
_{m}S_{m}(t)-u_2(t)S_{m}(t), \\ 
E_{m}^{\prime}(t)=B\beta _{hm}\frac{I(t)}{N}S_{m}(t)
-(\eta _{m}+\mu _{m})E_{m}(t)-u_2(t) E_{m}(t), \\ 
I_{m}^{\prime}(t)=\eta _{m}E_{m}(t)-\mu _{m}I_{m}(t)-u_2(t) I_{m}(t) \, .
\end{cases}
\label{eq4}
\end{equation}
The set $U$ of admissible control functions is defined by
\begin{equation*}
U=\left\{ \left. \left( u_{1}(\cdot),{u}_{2}(\cdot)\right) 
\in \left( L^\infty(0, t_f)\right)^2 \right\vert \text{ }0\leq
u_{1}\left( t\right),{u}_{2}\left( t\right) \leq u_{\max} \leq 0.5,\text{ }
\left[ 0,t_{f}\right] \right\}.
\end{equation*}

Our objective is twofold: first, to minimize the incidence of Zika virus infection 
among pregnant women, thereby reducing the risk of microcephaly in newborns; second, 
to decrease the mosquito population responsible for transmission. Achieving these goals 
involves implementing preventative measures while optimizing budget allocation 
to minimize costs. To achieve this, we consider the objective functional of the 
optimal control problem as follows:
\begin{equation}
J\left( I(\cdot),N_{m}(\cdot),u_{1,}(\cdot), u_{2}(\cdot)\right)
=\int_0^{t_f}\left( w_{1}I\left( t\right) 
+w_{2}N_{m}\left( t\right)+w_{3}u_{1}^{2}\left( t\right) 
+w_{4}u_{2}^{2}\left( t\right) \right) dt,
\label{eq5}
\end{equation}
{where the weight coefficients $w_{1}$ and $w_{2}$ are the weights for infected 
pregnant women and the mosquito population, respectively. Also, the coefficients 
$w_{3}$ and $w_{4}$ are measures of the cost of preventive interventions related 
to the controls $u_{1}$ and ${u}_{2}$, respectively.}

The optimal control problem consists of determining 
$$
X^{\ast }=\left(
S^{\ast },I^{\ast },W^{\ast },M^{\ast },A_{m}^{\ast },S_{m}^{\ast},
E_{m}^{\ast },I_{m}^{\ast }\right)
$$ 
associated with an admissible control
pair $u^{\ast }=\left( u_{1}^{\ast },{u}_{2}^{\ast }\right) \in U$ on the
time interval $\left[ 0,t_{f}\right] $, satisfying \eqref{eq4}, given initial
conditions $S\left( 0\right) ,$ $I\left( 0\right) ,$ $W\left( 0\right)$, 
$M\left( 0\right) ,$ $A_{m}\left( 0\right) ,$ $S_{m}\left( 0\right)$, 
$E_{m}\left( 0\right)$, and $I_{m}\left( 0\right) $ and minimizing the
objective functional \eqref{eq5}, i.e.,
\begin{equation}
\label{eq:minJ}
J\left( X^{\ast },u^{\ast }\right) 
=\underset{\left( X,u\right) 
\in \mathcal{X}\times U}{\min }J\left( X,u\right),
\end{equation}
where $\mathcal{X}$ is the set of admissible trajectories.
{We} achieve the necessary optimality conditions with the help of Pontryagin's
maximum principle (PMP) \cite{Pontryagin1987}. Note that the existence of optimal
controls is ensured by the convexity of the integrand of the objective
functional \eqref{eq5} with regard to the control functions 
$\left( u_{1},{u}_{2}\right)$ and the fact that the control system 
\eqref{eq4} satisfies a Lipschitz condition with respect 
to the state variables $S,$ $I,$ $W,$ $M,$ 
$A_{m},$ $S_{m},$ $E_{m}$, and $I_{m}$ 
\cite{Fleming1975,Cesari1983}.

\begin{Theorem}
Let $\left( u_{1}^{\ast },{u}_{2}^{\ast }\right) \in U$ be the optimal
controls that minimize the objective functional \eqref{eq5} and $\left(
S^{\ast },I^{\ast },W^{\ast },M^{\ast },A_{m}^{\ast },S_{m}^{\ast
},E_{m}^{\ast },I_{m}^{\ast }\right) $ be the optimal solution for the
dynamical system \eqref{eq4}. Thus, there are costate variables $\left(
\lambda _{1},\lambda _{2},\lambda _{3},\lambda _{4},\lambda _{5},\lambda
_{6},\lambda _{7},\lambda _{8}\right) $ {satisfying}
\begin{equation}
\small
\left\{ 
\begin{array}{lll}
\lambda _{1}^{\prime }(t) & = & \left( \lambda _{1}^{\prime }(t)-\lambda
_{2}^{\prime }(t)\right) \left( 1-u_{1}^{\ast }(t)\right) {\phi B\beta _{mh}
\frac{I_{m}^{\ast \prime }(t)}{\left( N^{\ast \prime }(t)\right) ^{2}}}
\left( I^{\ast \prime }(t)+W^{\ast \prime }(t)+M^{\ast \prime }(t)\right) \\ 
& + & {\left( \lambda _{1}^{\prime }(t)-\lambda _{3}(t)\right) (1-\phi )\tau
_{1}+\lambda _{1}(t)\mu _{h}},{-}\left( \lambda _{6}(t)-\lambda
_{7}(t)\right) {B\beta _{hm}\frac{I(t)}{\left( N^{\ast }(t)\right) ^{2}}
S_{m}^{\ast }}(t) \\ 
\lambda _{2}^{\prime }(t) & = & -w_{1}+\left( \lambda _{2}(t)-\lambda
_{1}(t)\right) \left( 1-u_{1}^{\ast }(t)\right) {\phi B\beta _{mh}
\frac{I_{m}^{\ast }(t)}{\left( N^{\ast }(t)\right) ^{2}}}S^{\ast }(t)+\lambda
_{2}(t){(\tau _{2}+\mu _{h})} \\ 
& - & {\lambda _{3}(t)(1-\psi )\tau _{2}-}\lambda _{4}(t){\psi \tau _{2}+}
\left( \lambda _{6}(t)-\lambda _{7}(t)\right) {B\beta _{hm}\frac{\left(
S^{\ast }(t)+W^{\ast }(t)+M^{\ast }(t)\right) }{\left( N^{\ast }(t)
\right)^{2}}S_{m}^{\ast }}(t) \\ 
\lambda _{3}^{\prime }(t) & = & \left( \lambda _{2}(t)-\lambda
_{1}(t)\right) \left( 1-u_{1}^{\ast }(t)\right) {\phi B\beta _{mh}
\frac{I_{m}^{\ast }(t)}{\left( N^{\ast }(t)\right) ^{2}}}S^{\ast }(t)
+{\lambda_{3}(t)\mu _{h}} \\ 
& - & \left( \lambda _{6}(t)-\lambda _{7}(t)\right) {B\beta _{hm}
\frac{I^{\ast }(t)}{\left( N^{\ast }(t)\right) ^{2}}S_{m}^{\ast }}(t) \\ 
\lambda _{4}^{\prime }(t) & = & \left( \lambda _{2}(t)-\lambda
_{1}(t)\right) \left( 1-u_{1}^{\ast }(t)\right) {\phi B\beta _{mh}
\frac{I_{m}^{\ast }(t)}{\left( N^{\ast }(t)\right) ^{2}}}S^{\ast }(t)
+{\lambda_{4}(t)\mu _{h}} \\ 
& - & \left( \lambda _{6}(t)-\lambda _{7}(t)\right) {B\beta_{hm}
\frac{I}{\left( N^{\ast }(t)\right) ^{2}}S_{m}^{\ast }}(t) \\ 
\lambda _{5}^{\prime }(t) & = & -\lambda _{5}(t)\left( {\mu _{b}\frac{1}{K}
\left( S_{m}^{\ast }(t)+E_{m}^{\ast }(t)+I_{m}^{\ast }(t)\right) +\left( \mu
_{A}+\eta _{A}\right) }\right) -\lambda _{6}(t){\eta _{A}} \\ 
\lambda _{6}^{\prime }(t) & = & -w_{2}-\lambda _{5}(t){\mu _{b}\left( 1
-\frac{A_{m}^{\ast }(t)}{K}\right) +}\left( \lambda _{6}(t)-\lambda
_{7}(t)\right) {B\beta _{hm}\frac{I^{\ast }(t)}{N^{\ast }(t)}+}\lambda
_{6}(t)\left( {\mu _{m}+u}_{2}^{\ast }(t)\right) \\ 
\lambda _{7}^{\prime }(t) & = & -w_{2}-\lambda _{5}(t){\mu _{b}\left( 1
-\frac{A_{m}^{\ast }(t)}{K}\right) +}\lambda _{7}(t)\left( {\eta _{m}
+\mu_{m}+u}_{2}^{\ast }(t)\right) -\lambda _{7}(t){\eta _{m}} \\ 
\lambda _{8}^{\prime }(t) & = & -w_{2}+\left( \lambda _{1}(t)-\lambda
_{2}(t)\right) \left( 1-u_{1}^{\ast }(t)\right) {\phi B\beta _{mh}
\frac{S^{\ast }(t)}{N^{\ast }(t)}-}\lambda _{5}(t){\mu _{b}\left( 1
-\frac{A_{m}^{\ast }(t)}{K}\right) } \\ 
& + & \lambda _{8}(t)\left( {\mu _{m}+u}_{2}^{\ast }(t)\right)
\end{array}
\right.  \label{eq6}
\end{equation}
with transversality conditions 
\begin{equation}
\lambda _{i}\left( t_{f}\right) =0,\text{ }\ i=1,2,\ldots,8.  \label{eq7}
\end{equation}
{In addition, }
\begin{equation}
\label{eq8}
\begin{cases}
u_{1}^{\ast }\left( t\right) 
= \min \left( \max \left( \frac{\left(
\lambda _{1}\left( t\right)-\lambda _{2}\left( t\right)\right) 
{\phi B\beta _{mh}\frac{I_{m}^{\ast }\left( t\right)}{N^{\ast }
\left( t\right)}}S^{\ast }\left( t\right)}{w_{3}},0\right) , u_{\max}\right), \\ 
u_{2}^{\ast }\left( t\right) 
= \min \left( \max \left( -\frac{\lambda
_{6}S_{m}^{\ast }\left( t\right)+\lambda _{7}I_{m}^{\ast }\left( t\right)
+\lambda _{8}E_{m}^{\ast }\left( t\right)}{w_{4}},0\right), u_{\max}\right).
\end{cases}
\end{equation}
\end{Theorem}

\begin{proof}
Using PMP \cite{Pontryagin1987}, we obtain
the necessary optimality conditions \eqref{eq6}--\eqref{eq8} that an optimal
solution must provide. We introduce the Hamiltonian $\mathcal{H}$ 
to form the necessary optimality conditions:

\begin{equation}
\small
\begin{array}{lll}
\mathcal{H}\left( t,X,u,\lambda\right) & = & w_{1}I+w_{2}N_{m}
+w_{3}u_{1}^{2}+w_{4}u_{2}^{2} \\ 
& + & \lambda _{1} \left( {\Lambda -}\left( 1-u_{1}\right) 
{\phi B\beta _{mh}\frac{I_{m}}{N}S
-(1-\phi )\tau _{1}S-\mu _{h}S}\right)  \\ 
& + & \lambda _{2} \left( {\left( 1-u_{1}\right) \phi B\beta
_{mh}\frac{I_{m}}{N}S
-(\tau _{2}+\mu _{h})I}\right)  \\ 
& + & \lambda _{3} \left( {(1-\phi )\tau _{1}S+(1-\psi )\tau
_{2}I-\mu _{h}W}\right)  \\ 
& +& \lambda _{4} \left( {\psi \tau _{2}I
-\mu _{h}M}\right) 
\\ 
&+ & \lambda _{5} \left( {\mu _{b}\left( 1-\frac{A_{m}}{K}
\right) \left( S_{m}+E_{m}+I_{m}\right) 
-\left( \mu _{A}+\eta _{A}\right)
A_{m}}\right)  \\ 
& + & \lambda _{6} \left( {\eta _{A}A_{m}-B\beta_{hm}
\frac{I}{N}S_{m}+\mu_{m}S_{m}
-u}_{2}{S_{m}}\right) \\ 
& + & \lambda _{7} \left( {B\beta _{hm}
\frac{I}{N} S_{m}
-(\eta _{m}+\mu _{m})E_{m}-u}_{2}{E_{m}}\right) \\ 
& + & \lambda _{8} \left( {\eta _{m}E_{m}
-\mu _{m}I_{m}-u}_{2}{I_{m}}\right).
\end{array}
\label{eq9}
\end{equation} 
{We} consider the following equations based on the Hamiltonian to obtain
the necessary optimality conditions of the problem:
\begin{itemize}
\item State equations:
\begin{equation}
\label{eq10}
\begin{gathered}
{\frac{dS}{dt}=}\frac{\partial \mathcal{H}}{\partial \lambda _{1}}\text{, }
{\frac{dI}{dt}=}\frac{\partial \mathcal{H}}{\partial \lambda _{2}}\text{, }
{\frac{dW}{dt}{=}\frac{\partial \mathcal{H}}{\partial \lambda _{3}},}\text{ }
{\frac{dM}{dt}{=}\frac{\partial \mathcal{H}}{\partial \lambda _{4}},} \\ 
{\frac{dA_{m}}{dt}{=}\frac{\partial \mathcal{H}}{\partial \lambda _{5}},}
\text{ }{\frac{dS_{m}}{dt}{=}\frac{\partial \mathcal{H}}{\partial 
\lambda_{6}},}\text{ }{\frac{dE_{m}}{dt}{=}\frac{\partial \mathcal{H}}{\partial
\lambda_{7}},}\text{ }{\frac{dI_{m}}{dt}=}\frac{\partial \mathcal{H}}
{\partial \lambda_{8}};
\end{gathered}
\end{equation}
\item Adjoint equations:
{\begin{equation}
\label{eq11}
\begin{gathered}
{\frac{d\lambda_{1}}{dt}=}-\frac{\partial \mathcal{H}}{\partial S}\text{, }
{\frac{d\lambda_{2}}{dt}=}-\frac{\partial \mathcal{H}}{\partial I}\text{, }
{\frac{d\lambda_{3}}{dt}=}-\frac{\partial \mathcal{H}}{\partial W},\text{ }
{\frac{d\lambda_{4}}{dt}=}-\frac{\partial \mathcal{H}}{\partial M}, \\ 
{\frac{d\lambda_{5}}{dt}=}-\frac{\partial \mathcal{H}}{\partial A_{m}},
\text{ }{\frac{d\lambda_{6}}{dt}=}-\frac{\partial \mathcal{H}}{\partial S_{m}},
\text{ }{\frac{d\lambda_{7}}{dt}=}-\frac{\partial \mathcal{H}}{\partial E_{m}},
\text{ }{\frac{d\lambda_{8}}{dt}=}-\frac{\partial \mathcal{H}}{\partial I_{m}},
\end{gathered}
\end{equation}}	
subject to transversality conditions
$\lambda _{i}(t_{f}) =0$, $i=1,2,\ldots,8$;
\item Minimality condition:
\begin{equation}
\label{eq12}
\mathcal{H}\left( t,X^*(t),u^*(t),\lambda(t)\right) 
= \min_{v \in [0,u_{\max}] \times [0,u_{\max}]}
\mathcal{H}\left(t,X^*(t),v,\lambda(t)\right).
\end{equation}			
\end{itemize}

We enforce these conditions to the Hamiltonian and see that the state equations \eqref{eq10}
correspond to the dynamical system \eqref{eq4}; we obtain the costate system 
\eqref{eq6}--\eqref{eq7} from the costate equations \eqref{eq11} and transversality
conditions; while we obtain the control functions \eqref{eq8} from
the minimality condition \eqref{eq12} of the PMP 
\cite{Pontryagin1987}. The proof is complete.
\end{proof}


\section{Numerical Simulations: Case Study in Brazil}
\label{sec:numerical}

In this section, we solve numerically the suggested optimal control problem, defined by 
\eqref{eq4}--\eqref{eq:minJ}. For this, we use the fourth-order Runge--Kutta
method (see, e.g., \cite{Book:Lenhart} for details). We consider the real data 
publicly available at the WHO, considered in \cite{Zika:Ndairou:etal}, 
of the confirmed cases of Zika in Brazil between 4 February  2016 and 10  November 2016. 
According to \cite{Zika:Ndairou:etal}, we consider as initial values $S_{0}$ = 2,180,686 
($S_{0}$ is the number of newborns corresponding to the simulation period) and the
number of births in the period, $I_{0}=1$, $M_{0}=0$, and $W_{0}=0$ for the
human female populations, and $A_{m0}=S_{m0}=I_{m0}=1.0903 \times 10^{6}$, and 
$E_{m0}=6.5421 \times 10^{6}$ for the mosquitoes populations. We obtain numerical simulations 
for a final time $t_f=160$ (weeks). Also, the parameters of the model \eqref{eq4} are listed 
in Table~\ref{parameters} (see \cite{Zika:Ndairou:etal}), and the maximum value of the control 
$u_1$ is assumed to be $u_{\max} = 0.5$.

\begin{figure}[H]
\hspace{-13mm}\includegraphics[width=1.1\linewidth]{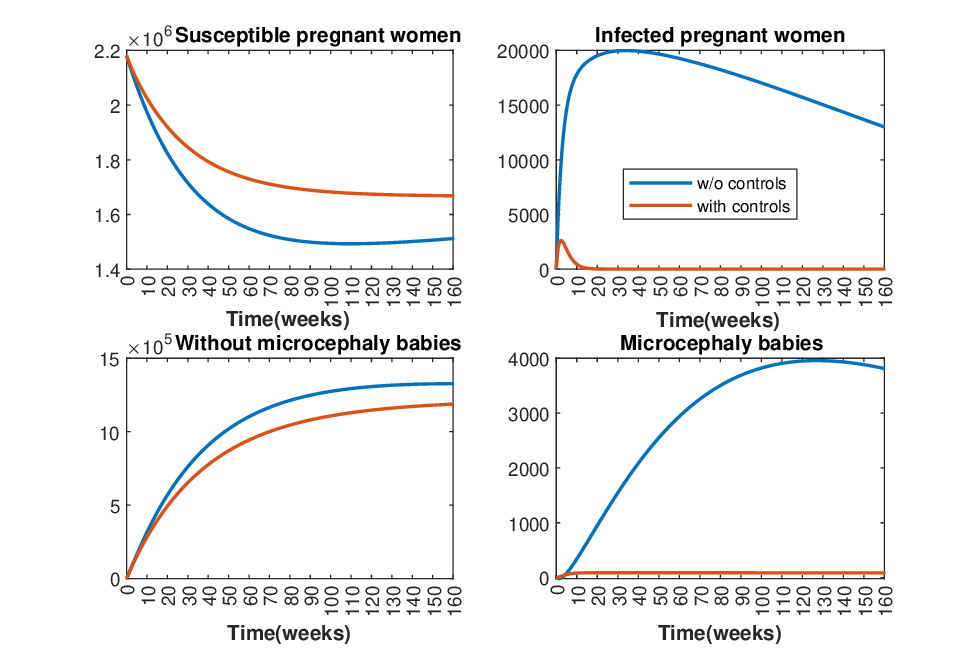}
\caption{Effect of control strategies on women population, with and without 
controls. (\textbf{{Top~ left}}): susceptible pregnant women $S$; {(\textbf{{Top right}}):} 
infected pregnant women $I$; (\textbf{{Bottom left}}): 
women who gave birth to babies without microcephaly $W$; 
(\textbf{{Bottom right}}): women who gave birth to babies with microcephaly $M$.} 
\label{fig1}
\end{figure}

\begin{figure}[H]
\hspace{-12mm}\includegraphics[width=1.1\linewidth]{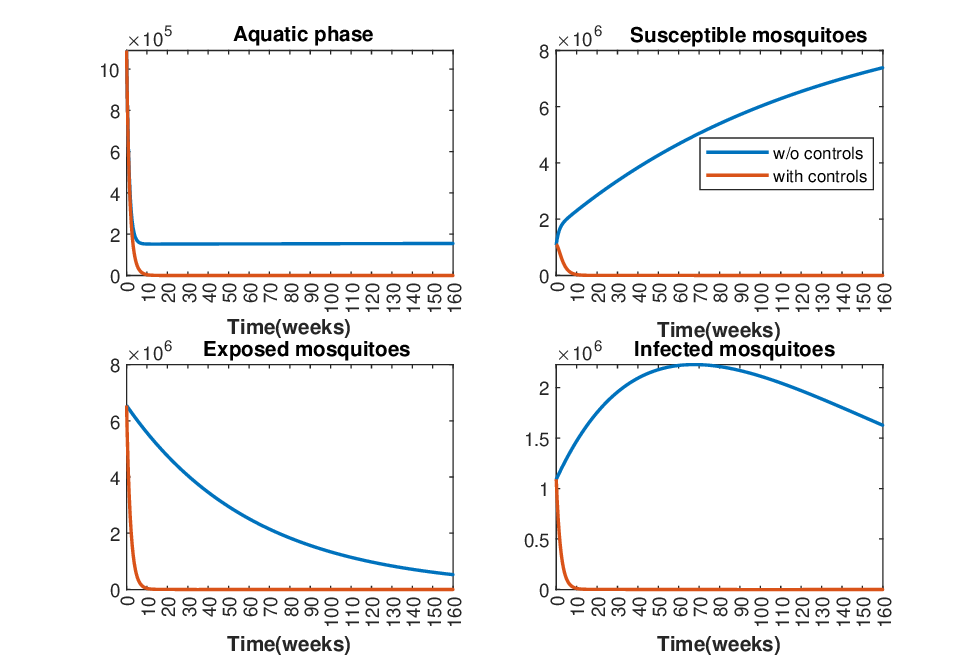} 
\caption{Effect of control strategies on mosquitoes, with and without controls.
(\textbf{Top left}): mosquitoes in the aquatic phase $A_m$; (\textbf{Top right}): susceptible mosquitoes 
$S_m$; (\textbf{Bottom left}): exposed mosquitoes $E_m$; (\textbf{Bottom right}): infected mosquitoes $I_m$.} 
\label{fig2}
\end{figure}

\begin{figure}[H]
\hspace{-12mm}\includegraphics[width=1\linewidth]{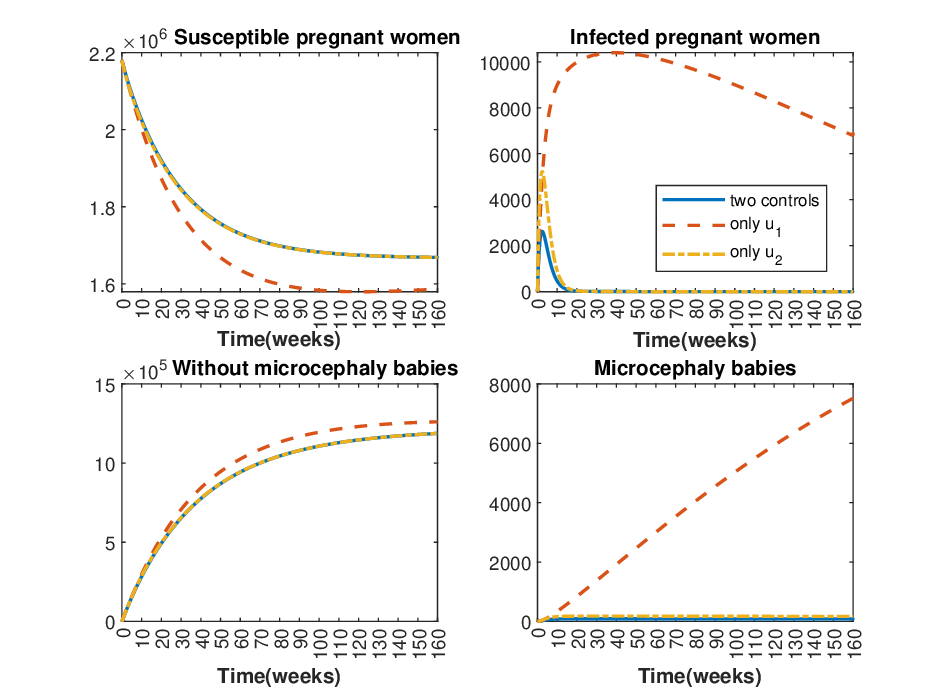} 
\caption{Comparative impact of the control strategies on women population, 
considering both controls and one single control $u_1$ or $u_2$. (\textbf{Top left}): susceptible 
pregnant women $S$; (\textbf{Top right}): infected pregnant women $I$; (\textbf{Bottom left}): women who gave 
birth to babies without microcephaly $W$; (\textbf{Bottom right}): women who gave birth to babies 
with microcephaly $M$.} 
\label{fig3}
\end{figure}

\begin{figure}[H]
\hspace{-12mm}\includegraphics[width=1\linewidth]{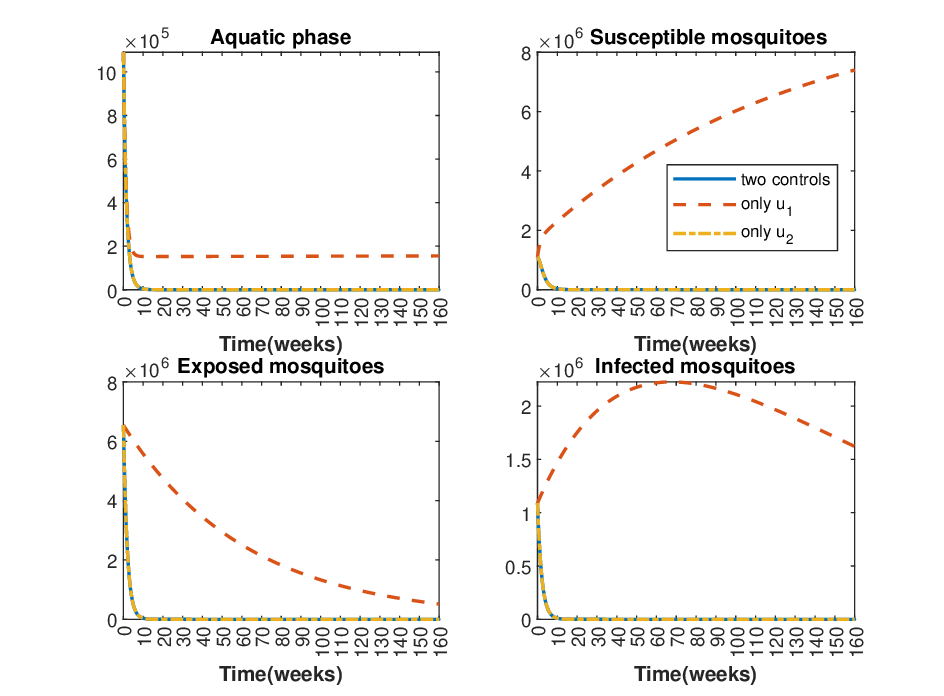} 
\caption{Comparative impact of the control strategies on mosquitoes, considering 
both controls and one single control $u_1$ or $u_2$. (\textbf{Top left}): mosquitoes in the 
aquatic phase $A_m$; (\textbf{Top right}): susceptible mosquitoes $S_m$; (\textbf{Bottom left}): exposed 
mosquitoes $E_m$; (\textbf{Bottom right}): infected mosquitoes $I_m$.} 
\label{fig4}
\end{figure}

\begin{figure}[H]
\hspace{-6mm}\includegraphics[width=1\linewidth]{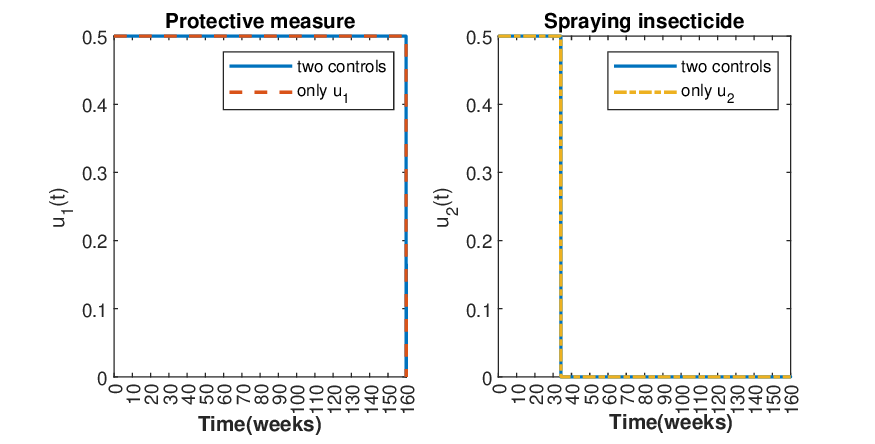} 
\caption{Control strategies $u_1$ and $u_2$.} 
\label{fig5}
\end{figure}

\begin{figure}[H]
\hspace{-12mm}\includegraphics[width=1.1\linewidth]{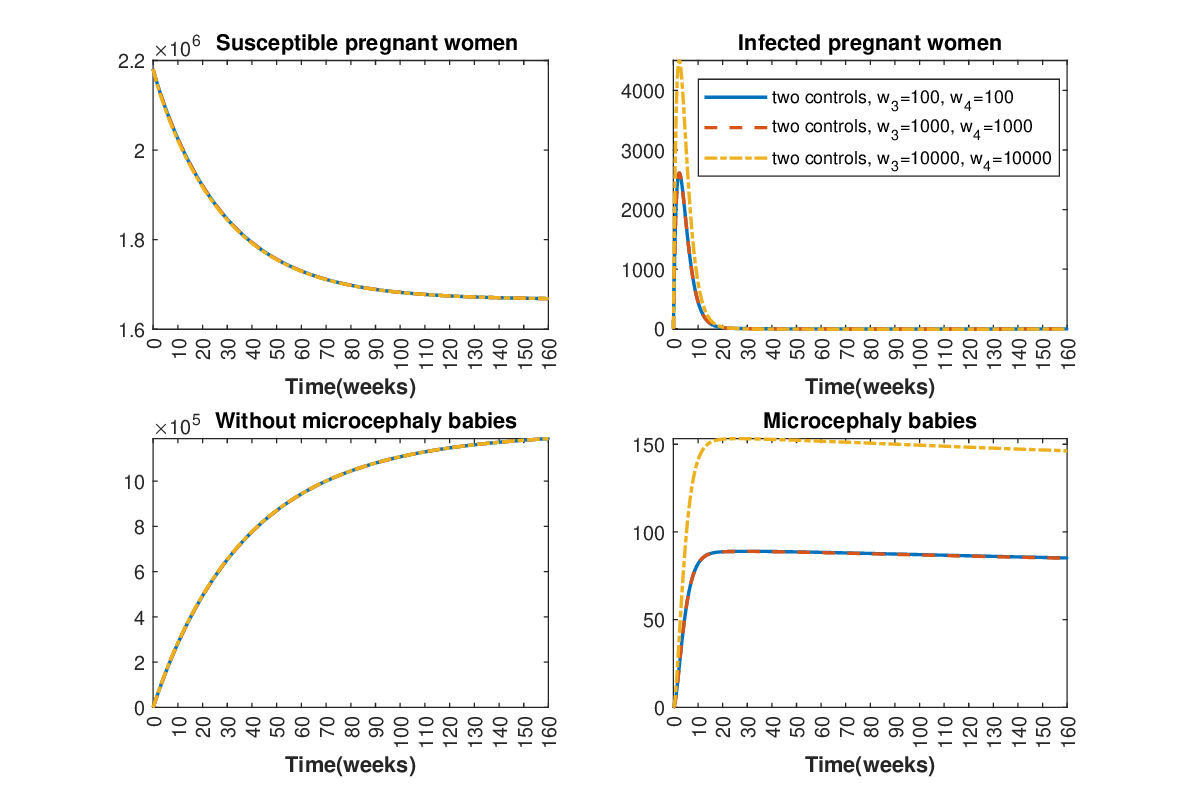} 
\caption{Comparative impact of control strategies with different weight coefficients 
on women population, considering two controls. (\textbf{Top left}): susceptible pregnant women $S$;
(\textbf{Top right}): infected pregnant women $I$; (\textbf{Bottom left}): women who gave birth to babies 
without microcephaly $W$; (\textbf{Bottom right}): women who gave birth to babies with microcephaly $M$.}
\label{fig6}
\end{figure}

\begin{figure}[H]
\hspace{-10mm}\includegraphics[width=1.1\linewidth]{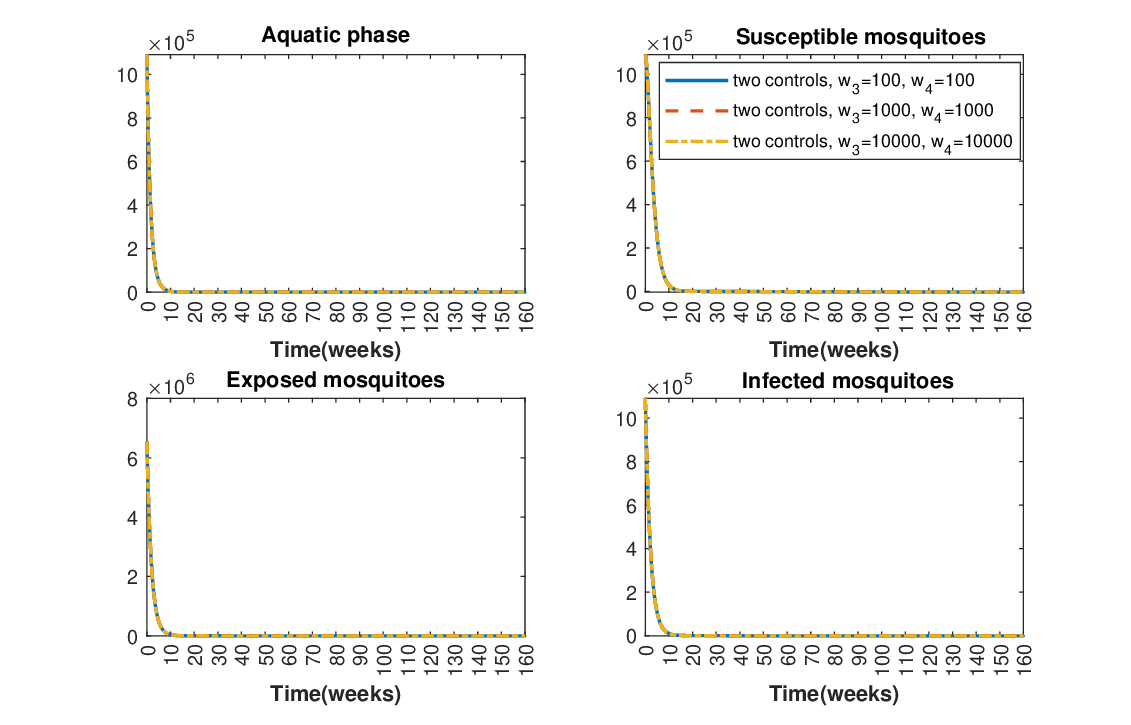} 
\caption{Comparative impact of control strategies with different weight coefficients 
on mosquitoes, considering two controls. (\textbf{Top left}): mosquitoes in the 
aquatic phase $A_m$; (\textbf{Top right}): susceptible mosquitoes $S_m$; 
(\textbf{Bottom left}): exposed mosquitoes $E_m$; 
(\textbf{Bottom right}): infected mosquitoes $I_m$.}
\label{fig7}
\end{figure}

\begin{figure}[H]
\hspace{-10mm}\includegraphics[width=1\linewidth]{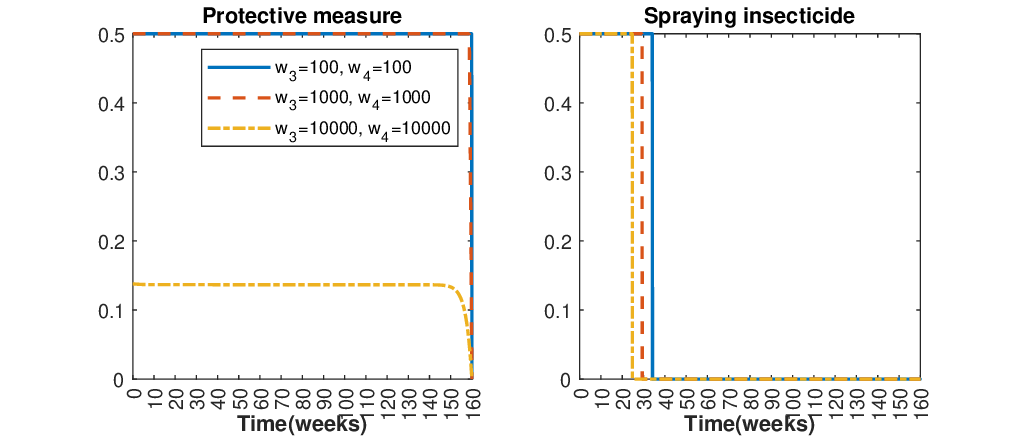} 
\caption{Control strategies $u_1$ and $u_2$ with weight coefficients $w_{1} = w_{2} = 10$
and $w_{3} = w_{4}$ with $w_{3}, w_{4} \in \{100, 1000, 10000\}$.} 
\label{fig8}
\end{figure}

\begin{figure}[H]
\hspace{-10mm}\includegraphics[width=1.0\linewidth]{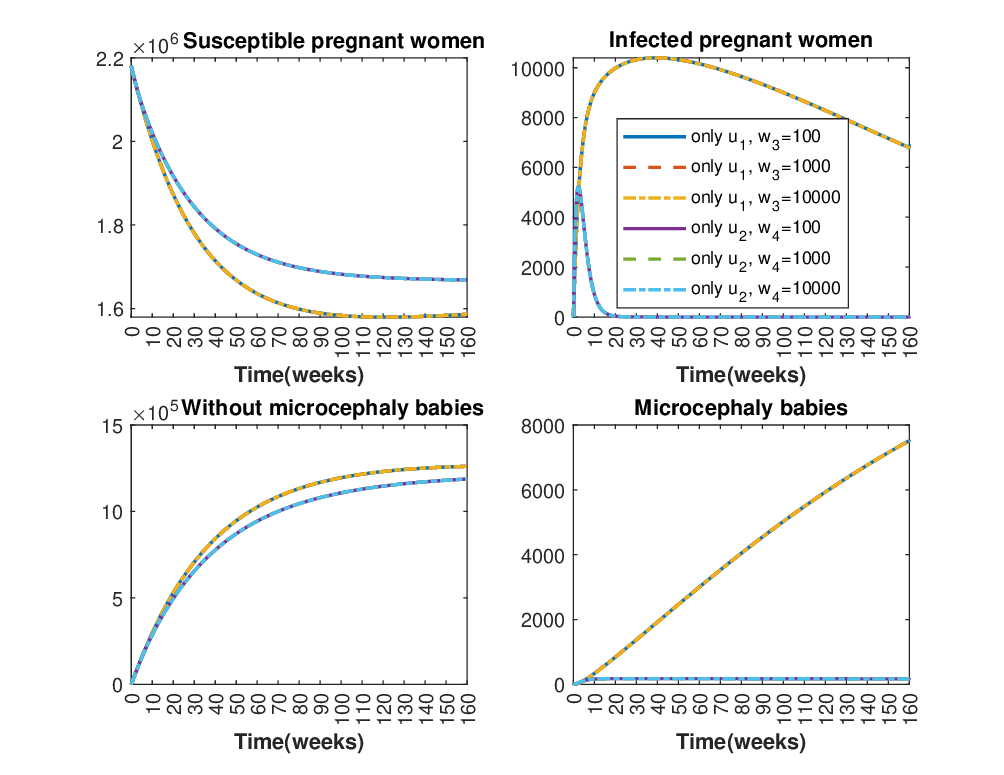} 
\caption{Comparative impact of control strategies with different weight coefficients 
on women population, considering one single control $u_1$ or $u_2$. 
(\textbf{Top left}): susceptible pregnant women $S$; 
(\textbf{Top right}): infected pregnant women $I$; 
(\textbf{Bottom left}): women who gave birth to babies without microcephaly $W$; 
(\textbf{Bottom right}): women who gave birth to babies with microcephaly $M$.}
\label{fig9}
\end{figure}

\begin{figure}[H]
\hspace{-10mm}\includegraphics[width=1.0\linewidth]{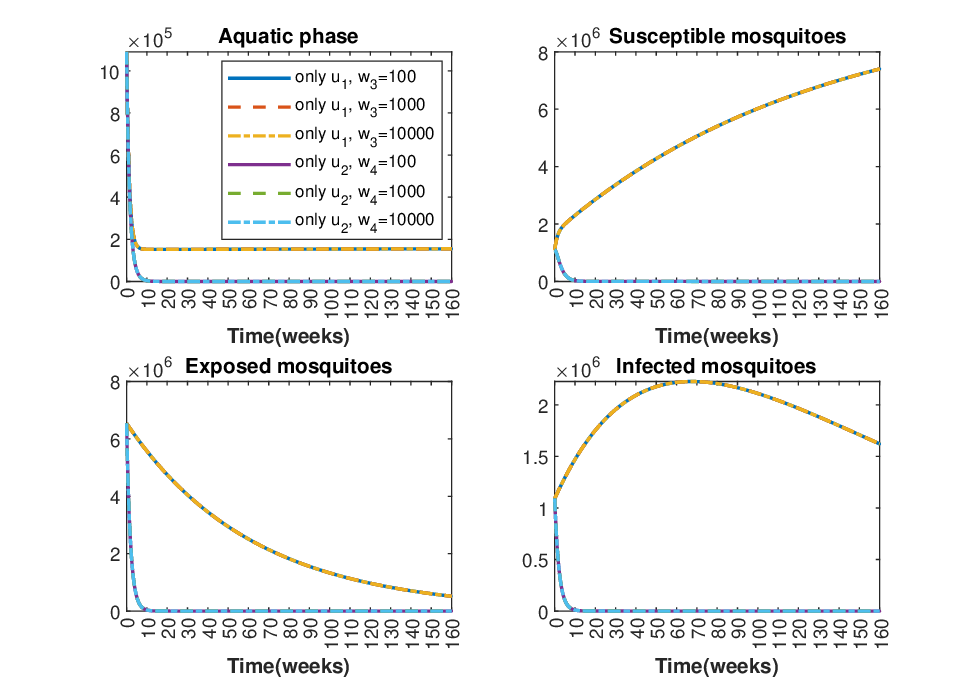} 
\caption{Comparative 
\textls[-25]{ impact of control strategies with different weight coefficients 
on mosquitoes, considering one single control $u_1$ or $u_2$. 
(\textbf{Top left}): mosquitoes in the aquatic phase $A_m$; 
(\textbf{Top right}): susceptible mosquitoes $S_m$; 
(\textbf{Bottom left}): exposed mosquitoes $E_m$; 
(\textbf{Bottom right}): infected \mbox{mosquitoes $I_m$.}}}
\label{fig10}
\end{figure}

\begin{figure}[H]
\hspace{-10mm}\includegraphics[width=1\linewidth]{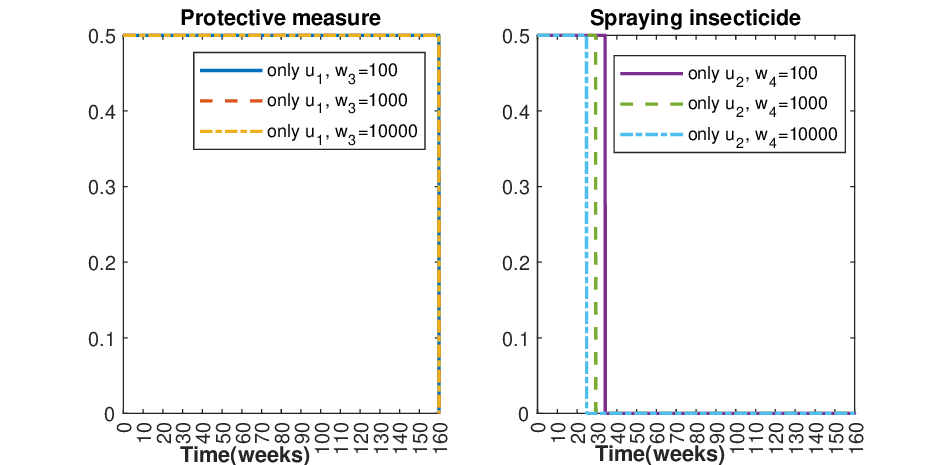} 
\caption{Control strategy $u_1$ with $w_{1} = w_{2} = 10$
and $w_{3} \in \{100, 1000, 10000\}$ ($u_2 \equiv 0$);
control strategy $u_2$ with $w_{1} = w_{2} = 10$
and $w_{4} \in \{100, 1000, 10000\}$ ($u_1 \equiv 0$).}
\label{fig11}
\end{figure}

\begin{table}[H]
\caption{{Parameter values for system 
\eqref{zikamodel1}--\eqref{zikamodel2} and weight coefficients for \eqref{eq5}.}}
\label{parameters}
\begin{tabular}{llll}
\toprule
{\footnotesize {\textbf{Symbol}}} & {\footnotesize {\textbf{Description}}} 
& {\footnotesize {\textbf{Value}}} & {{\footnotesize {\textbf{References}}}} \\ \midrule
{\footnotesize {$\Lambda$}} & {\footnotesize {new pregnant women (per week)}} & 
{\footnotesize {3,000,000/52 }} & {\footnotesize {\cite{Zika:Ndairou:etal}}} \\ 
{\footnotesize {$\phi$}} & {\footnotesize {fraction of $S$ that become infected}}
& {\footnotesize {$0.459$}} & {\footnotesize {\cite{Zika:Ndairou:etal}}}\\ 
{\footnotesize {$B$}} & {\footnotesize {average daily biting (per day)}} & 
{\footnotesize {$1$}} & {\footnotesize {\cite{Zika:Ndairou:etal}}} \\ 
{\footnotesize {$\beta_{mh}$}} & {\footnotesize {transmission probability from 
$I_m$ (per bite)}} & {\footnotesize {$0.6$}} & {\footnotesize {\cite{Zika:Ndairou:etal}}}\\ 
{\footnotesize {$\tau_{1}$}} & {\footnotesize {rate at which $S$ give birth (in
weeks)}} & {\footnotesize {$37$}} & {\footnotesize {\cite{Zika:Ndairou:etal}}}\\ 
{\footnotesize {$\tau_{2}$}} & {\footnotesize {rate at which $I$ give birth (in
weeks)}} & {\footnotesize {$1/25$}} & {\footnotesize {\cite{Zika:Ndairou:etal}}}\\ 
{\footnotesize {$\mu_{h}$}} & {\footnotesize {natural death rate}} & 
{\footnotesize {$1/50$}} & {\footnotesize {\cite{Zika:Ndairou:etal}}} \\ 
{\footnotesize {$\psi$}} & {\footnotesize {fraction of $I$ that gives birth to
babies with a neurological disorder}} 
& {\footnotesize {$0.133$}} & {\footnotesize {\cite{Zika:Ndairou:etal}}}\\ 
{\footnotesize {$\beta_{hm}$}} & {\footnotesize {transmission probability from 
$I_h$ (per bite)}} & {\footnotesize {$0.6$}}& {\footnotesize {\cite{Zika:Ndairou:etal}}} \\ 
{\footnotesize {$\mu_b$}} & {\footnotesize {number of eggs at each deposit per
capita (per day)}} & {\footnotesize {$80$}}& {\footnotesize {\cite{Zika:Ndairou:etal}}} \\ 
{\footnotesize {$\mu_A$}} & {\footnotesize {natural mortality rate of larvae
(per day)}} & {\footnotesize {$1/4$}} & {\footnotesize {\cite{Zika:Ndairou:etal}}} \\ 
{\footnotesize {$\eta_A$}} & {\footnotesize {maturation rate from larvae to
adult (per day)}} & {\footnotesize {$0.5$}}& {\footnotesize {\cite{Zika:Ndairou:etal}}} \\ 
{\footnotesize {$1/\eta_m$}} & {\footnotesize {extrinsic incubation period (in
days)}} & {\footnotesize {$125$}} & {\footnotesize {\cite{Zika:Ndairou:etal}}} \\ 
{\footnotesize {$1/\mu_m$}} & {\footnotesize {average lifespan of adult
mosquitoes (in days)}} & {\footnotesize {$125$}} & {\footnotesize {\cite{Zika:Ndairou:etal}}} \\ 
{\footnotesize {$K$}} & {\footnotesize {maximal capacity of larvae}} & 
{\footnotesize {$1.09034 \times 10^{6}$}} & {\footnotesize {\cite{Zika:Ndairou:etal,Rodrigues}}} \\ 
{\footnotesize {$w_{1}$}} 
& {{\footnotesize {the weight coefficient for infected pregnant women}}} & 
{{\footnotesize {$10$ }}} 
& {{\footnotesize {assummed}}}\\ 
{{\footnotesize {$w_{2}$}}} 
& {{\footnotesize {the weight coefficient for the mosquito population}}} & 
{{\footnotesize {$10$ }}} 
& {{\footnotesize {assummed}}} \\ 
{{\footnotesize {$w_{3}$}}} 
& {{\footnotesize {the weight coefficient for the cost of protective measure}}} & 
{{\footnotesize {$100$ }}} 
& {{\footnotesize {assummed}}} \\ 
{{\footnotesize {$w_{4}$}}} 
& {{\footnotesize {the weight coefficient for the cost of spraying insecticide}}} & 
{{\footnotesize {$100$ }}} 
& {{\footnotesize {assummed}}} \\ \bottomrule
\end{tabular}
\end{table}
In Figures~\ref{fig1} and \ref{fig2}, graphical representations depict 
the transmission dynamics of Zika within women and mosquito populations, 
with and without implementation of control measures. These visualizations 
illustrate the efficacy of control measures in averting cases of microcephaly 
by diminishing the count of infected pregnant women and eradicating the mosquito 
population responsible for Zika transmission. In Figures \ref{fig3} and \ref{fig4}, 
we observe that employing the suggested combination of control measures robustly 
hampers Zika transmission. {The mosquito population is essentially eradicated by the 
10th week through the implementation of insecticide spraying control measures. Furthermore, 
the combined effect of both control measures resulted in a notable reduction in 
infection cases in pregnant women by the 10th week. A comparison of the control strategies 
reveals that the two control strategies are the most effective. However, the results also 
demonstrate that the spraying insecticide measure has a notable impact.}
Lastly, in Figure~\ref{fig5}, while control $u_{1}$ remains in effect during all time 
window, control $u_{2}$ ceases to be effective around the 40th week due to its near-complete 
elimination of the mosquito population. {Additionally, Figures~\ref{fig6}--\ref{fig11} 
show comparative results of control strategies, considering different weight coefficients.
As illustrated in Figures~\ref{fig6} and \ref{fig7}, an increase in the weight coefficients 
of the two control strategies is associated with a rise in the number of babies born with 
microcephaly. This phenomenon occurs concurrently with no discernible impact on the mosquito 
population. As shown in Figure~\ref{fig8}, this phenomenon can be attributed to the diminished 
rate of exertion associated with the protective measure control in response to an elevated 
weight coefficient. Although a similar scenario would be expected for only one control strategy, 
it affects only the behavior of the spraying insecticide control rather than the behavior of the 
controlled system (see Figures~\ref{fig9}--\ref{fig11}). As a result, it reveals that in all cases, 
the two control strategies are the most robust for vertical Zika transmission.}

\vspace{-20PT}

\vspace{-20PT}
\vspace{-10PT}
\vspace{-15PT}
\vspace{-10PT}
\vspace{-20PT}


\section{Conclusions}
\label{sec:conclusions}

In this study, we have proposed control strategies aimed at thwarting 
the transmission of Zika virus, known for  precipitating neurological 
disorders like microcephaly. {For this purpose, an optimal control 
problem has been formulated for a model representing the vertical 
transmission of the Zika virus from infected mothers to infants in Brazil. 
Our primary objective has been to curtail the incidence of infection among 
pregnant women and diminish the mosquito population, all while ensuring 
cost-effectiveness in implementing these strategies. The optimal control 
problem has been solved via Pontryagin's maximum principle. Finally, the 
numerical results have been obtained using the fourth-order \mbox{Runge--Kutta }
method with the help of the MATLAB (2021b) numeric computing environment.} 
Through numerical simulations, we have demonstrated that the adoption 
of these control measures has led to a consistent reduction in the count 
of infected pregnant women from the outset of the intervention, 
consequently resulting in a decline in cases of microcephaly.


\vspace{6pt} 

\authorcontributions{Conceptualization, C.J.S. and D.F.M.T.; 
methodology, D.Y., C.J.S., and D.F.M.T.; 
software, D.Y.; 
validation, C.J.S. and D.F.M.T.; 
formal analysis, D.Y., C.J.S., and D.F.M.T.; 
investigation, D.Y., C.J.S., and D.F.M.T.; 
writing---original draft preparation, D.Y., C.J.S., and D.F.M.T.; 
\mbox{writing---review} and editing, D.Y., C.J.S., and D.F.M.T.; 
visualization, D.Y.; 
supervision, D.F.M.T. 
All authors have read and agreed to the published version of the manuscript.}

\funding{This research was funded by 
Funda\c{c}\~{a}o para a Ci\^{e}ncia e a Tecnologia (FCT),
grant numbers: UIDB/04106/2020 (\url{https://doi.org/10.54499/UIDB/04106/2020}) 
and UIDP/04106/2020 (\url{https://doi.org/10.54499/UIDP/04106/2020});	
and project Mathematical Modelling of Multiscale
Control Systems: Applications to Human Diseases (CoSysM3), 
reference 2022.03091.PTDC (\url{https://doi.org/10.54499/2022.03091.PTDC}).}

\dataavailability{Data are contained within the article.} 

\acknowledgments{\textls[-25]{{D.Y. is grateful to the 2211/A General
Domestic Doctoral Scholarship supported by TUBITAK and an Erasmus+ Grant}}. 
The authors would like to sincerely thank the two referees for contributing 
to the initially submitted version with several constructive questions, 
comments, and suggestions.}

\conflictsofinterest{The authors declare no conflicts of interest.} 


\newpage

\begin{adjustwidth}{-\extralength}{0cm}
	
\reftitle{References}



\PublishersNote{}

\end{adjustwidth}



\begin{thebibliography}{999}
	
\bibitem{WHO1} 
{Kindhauser, M.K.; Allen, T.; Frank, V.; Santhana, R.S.; Dye, C.}
Zika: The origin and spread of a mosquito-borne virus. 
\emph{Bull. World Health Organ.} \textbf{2016},   \emph{94}, {675--686.}

\bibitem{WHO2} 
World Health Organization. 
World Health Organization: Fact Sheets About Vector-Borne Diseases.
Available online: \url{https://www.who.int/news-room/feature-stories/detail/the-history-of-zika-virus}  
(accessed on November 5, 2024).

\bibitem{WHO3} 
World Health Organization. 
Zika Virus, Key Facts. 
Available online:  \url{https://www.who.int/news-room/fact-sheets/detail/zika-virus} 
(accessed on November 5, 2024). 

\bibitem{Martins} 
Almeida, R.; Brito da Cruz, A.M.; Martins, N.; Monteiro, M.T.T. 
An epidemiological MSEIR model described by the Caputo fractional derivative. \emph{Int. 
J. Dyn. Control} \textbf{2019}, \emph{7}, 776--784.

\bibitem{Cholera} 
Lemos-Pai\~{a}o, A.P.; Maurer, H.; Silva, C.J.; Torres, D.F.M. 
A SIQRB delayed model for cholera and optimal control treatment. 
{\emph{Math. Model. Nat. Phenom. }\textbf{2022}, \emph{17}, 25.}
{\tt arXiv:2206.12688}

\bibitem{Monkeypox} 
Yap\i \c{s}kan, D.; Yurto\u{g}lu, M.; Avc\i\,  D.; Ero\u{g}lu, B.B.\.{I}; Bonyah, E. 
A Novel Model for Monkeypox Disease: System Analysis and Optimal Preventive Strategies. 
\emph{Iran.  J. Sci.} \textbf{2023}, \emph{5}, 1--13.

\bibitem{Agusto} 
Agusto, F.B.; Bewick, S.; Fagan, W.F.
Mathematical model of Zika virus with vertical transmission. 
\emph{Infect. Dis. Model.} \textbf{2017}, \emph{2}, 244--267.

\bibitem{Shah} 
Shah, N.H.; Patel, Z.A.; Yeolekar, B.M. 
Preventions and controls on congenital transmissions of Zika: Mathematical analysis.
\emph{Appl.  Math.} \textbf{2017}, \emph{8}, {500--519.}

\bibitem{Zika:Ndairou:etal} 
{Nda\"{\i}rou, F.; Area, I.; Nieto, J.J.; Silva, C.J.; Torres, D.F.M.}  
Mathematical modeling of Zika disease in pregnant women and newborns with microcephaly in Brazil. 
\emph{Math. Meth. Appl. Sci}. \textbf{2018}, \emph{41},  8929--8941.
{\tt arXiv:1711.05630}

\bibitem{Mosquito} 
Ghaffari, P.; Silva, C.J.; Torres, D.F.M.
Mathematical Models and Optimal Control in Mosquito Transmitted Diseases.
In \emph{Bio-mathematics, Statistics, and Nano-Technologies: Mosquito Control Strategies}; 
Ghaffari, P., Ed.; Chapman and Hall/CRC: Boca Raton, FL, USA, 2023; pp.~143--156.

\bibitem{COVID-19} 
Ero\u{g}lu, B.B.\.{I}; Yap\i\c{s}kan, D. 
Optimal Strategies to Prevent COVID-19 from Becoming a Pandemic. 
In \emph{Mathematical Modeling and Intelligent Control for Combating Pandemics};  
Hammouch, Z., Lahby, M., Baleanu, D., Eds.;  
{Springer Optimization and Its Applications};  
Springer: Cham, Switzerland, 2023; pp. 39--55.

\bibitem{cancer} 
Esmaili, S.; Eslahchi, M.R.; Torres, D.F.M. 
Optimal control for a nonlinear stochastic PDE model of cancer growth.
{\emph{Optimization} \textbf{2024}, \emph{73},  2745--2789.}
{\tt arXiv:2307.09574}

\bibitem{epidemicmodel} 
Yap\i\c{s}kan, D.; Ero\u{g}lu, B.B.\.{I}. 
Fractional optimal control of a generalized SIR epidemic model with vaccination 
and treatment. In  \emph{Fractional Dynamics in Natural Phenomena and Advanced Technologies}; 
Baleanu, D., Hristov, J., Eds.; 
Cambridge Scholars Publishing: {Cambridge, MA, USA}, 2024; pp. 131--150.

\bibitem{Wang} 
Wang, X.; Shen, M.; Xiao, Y.; Rong, L. 
Optimal control and cost-effectiveness analysis of a Zika virus infection model 
with comprehensive interventions. 
\emph{Appl.  Math.  Comput.} \textbf{2019}, \emph{359}, 165--185.

\bibitem{Okyere} 
Okyere, E.; Olaniyi, S.; Bonyah, E. 
Analysis of Zika virus dynamics with sexual transmission route 
using multiple optimal controls. 
\emph{Sci.  Afr.} \textbf{2020}, \emph{9}, e00532.

\bibitem{Ali} 
Ali, A.; Iqbal, Q.; Asamoah, J.K.K.; Islam, S. 
Mathematical modeling for the transmission potential 
of Zika virus with optimal control strategies. 
\emph{ Eur.  Phys. J. Plus} \textbf{2022}, \emph{137},  146.

\bibitem{CDC:protect:Zika}
Centers for Disease Control and Prevention.
Available online: \url{https://www.cdc.gov/zika/prevention/index.html} 
(accessed on November 5, 2024). 

\bibitem{Pontryagin1987} 
Pontryagin, L.; Boltyanskii, V.; Gramkrelidze,  R.; Mischenko, E.  
\emph{The Mathematical Theory of Optimal Processes}; 
Wiley Interscience: {Hoboken, NJ, USA, }1962.

\bibitem{Fleming1975} 
{Fleming, W.H.; Rishel, R.W.} 
\emph{Deterministic and Stochastic Optimal Contro}l; 
Springer: New York, NY, USA,  1975.

\bibitem{Cesari1983} 
{Cesari, L.} 
\emph{Optimization---Theory and Applications: Problems 
with Ordinary Differential Equations}; Springer: New York, NY, USA, 1983.

\bibitem{Book:Lenhart}
{Lenhart, S.; Workman, J.T.} 
\emph{Optimal Control Applied to Biological Models}; 
Chapman \& Hall/CRC: Boca Raton, FL, USA, 2007.

\bibitem{Rodrigues}
{Rodrigues, H.S.; Monteiro, M.T.T.; Torres, D.F.M.; Zinober, A. 
Dengue disease, basic reproduction number and control. 
\emph{Int. J. Comput. Math.} \textbf{2012}, \emph{89}, 334--346.}
{\tt arXiv:1103.1923}

\end{thebibliography}
\end{document}